\documentclass{amsart}
\usepackage[nohug]{diagrams}
\usepackage{amscd}
\usepackage[mathscr]{eucal}
\usepackage{amsfonts}
\usepackage{amsmath}
\usepackage{amsthm}
\usepackage{amssymb}
\usepackage{latexsym}
\usepackage{tabularx,array}
\newfont{\msam}{msam10}

\newtheorem{theorem}[]{Theorem}
\newtheorem{proposition}[]{Proposition}
\newtheorem{corollary}[]{Corollary}
\newtheorem{lemma}[]{Lemma}
\theoremstyle{definition}
\newtheorem{definition}[]{Definition}
\newtheorem{remark}[]{Remark}
\newtheorem{example}[]{Example}

\let\nc\newcommand

\def\bthm{\begin{theorem}}
\def\ethm{\end{theorem}}
\def\blemma{\begin{lemma}}
\def\elemma{\end{lemma}}
\def\bproof{\begin{proof}}
\def\eproof{\end{proof}}
\def\bprop{\begin{proposition}}
\def\eprop{\end{proposition}}
\def\bcor{\begin{corollary}}
\def\ecor{\end{corollary}}
\nc{\la}{\label}
\def\Z{\mathbb{Z}}

\def\M{\mathbb{M}}
\def\X{\mathbb{X}}

\def\Alg{\mathtt{Alg}}
\def\Mod{\mathtt{Mod}}
\def\GrMod{\mathtt{GrMod}}
\def\Bimod{\mathtt{Bimod}}
\def\cAlg{\mathtt{Comm\,Alg}}
\def\gcAlg{\mathtt{GrComm\,Alg}}
\def\Sets{\mathtt{Sets}}
\def\DGA{\mathtt{DGA}}
\def\cDGA{\mathtt{Comm\,DGA}}
\def\DGMod{\mathtt{DG\,Mod}}
\def\DGBimod{\mathtt{DG\,Bimod}}
\def\D{{\mathsf D}}
\def\rtv#1{\!\!\sqrt[V]{#1}}
\nc{\Hom}{{\rm{Hom}}}
\nc{\HOM}{\underline{\rm{Hom}}}
\nc{\DER}{\underline{\rm{Der}}}
\nc{\Ext}{{\rm{Ext}}}
\nc{\Rep}{{\rm{Rep}}}
\nc{\DRep}{{\rm{DRep}}}
\nc{\NCRep}{\widetilde{\rm{Rep}}}
\nc{\RAct}{{\rm{RAct}}}
\nc{\EXT}{\underline{\rm{Ext}}}
\nc{\TOR}{\underline{\rm{Tor}}}
\def\H{\mathrm H}
\def\T{\mathrm T}
\nc{\End}{{\rm{End}}}
\nc{\GL}{{\rm{GL}}}
\nc{\PGL}{{\rm{PGL}}}
\nc{\SL}{{\rm{SL}}}
\nc{\PSL}{{\rm{PSL}}}
\nc{\ad}{{\rm{ad}}}
\nc{\dlim}{\varinjlim}
\def\eA{A^{\mbox{\scriptsize{\rm{e}}}}}

\def\eR{R^{\mbox{\scriptsize{\rm{e}}}}}

\newcommand{\Spec}{{\rm{Spec}}}

\newcommand{\Der}{{\rm{Der}}}

\newcommand{\into}{\,\,\hookrightarrow\,\,}

\begin{document}

\title{A Simple Construction of Derived Representation Schemes}
\author{Yuri Berest}
\address{Department of Mathematics,
 Cornell University, Ithaca, NY 14853-4201, USA}
\email{berest@math.cornell.edu}
\author{George Khachatryan}
\address{Department of Mathematics,
Cornell University, Ithaca, NY 14853-4201, USA}
\email{georgek@math.cornell.edu}
\author{Ajay Ramadoss}
\address{Departement Mathematik,
Eidgenossische TH Z\"urich,
8092 Z\"urich, Switzerland}
\email{ajay.ramadoss@math.ethz.ch}
\begin{abstract}
We present a simple algebraic construction of the (non-abelian) derived functors $ \DRep_n^\bullet(A) $
of the representation scheme $ \Rep_n(A) $, parametrizing the $n$-dimensional representations
of an associative algebra $A$. We construct a related derived version of the representation functor
introduced recently by M.~Van den Bergh \cite{vdB} and,  as an application, compute the derived tangent spaces
$ \T_\varrho^{\bullet} \DRep_n(A) $ to $ \Rep_n(A) $. We prove that our construction of $ \DRep_n^\bullet(A) $
agrees with an earlier construction of derived action spaces, due to I.~Ciocan-Fontanine and M.~Kapranov
\cite{CK}; however, our approach, proofs and motivation are quite different.
This paper is mainly a research announcement; detailed proofs and applications will appear elsewhere.
\end{abstract}
\maketitle

\section{Introduction}
\la{S1}
Let $A$ be an associative unital algebra over a field $k$. The classical representation scheme,
parametrizing the $n$-dimensional representations of $A$, can be defined as
the functor on the category of commutative algebras
\begin{equation}
\la{rep}
\Rep_n(A):\ \cAlg_{k} \to \Sets\ ,\quad B \mapsto \Hom_{\Alg}(A,\, \M(n,B))\ ,
\end{equation}
where $ \M(n,B) $ denotes the ring of $\,n \times n \,$ matrices over the commutative
algebra $B$. A natural way to prove representability of this functor goes back to the
work of Bergman \cite{B} and Cohn \cite{C}: the idea is to extend \eqref{rep}
from $\,\cAlg_k\,$ to the category of {\it all} associative $k$-algebras:
\[
\begin{diagram}[small, tight]
\cAlg_{k} &  \rTo^{\ \Rep_n(A)\quad }        & \Sets \\
\dInto^{} &  \ruDotsto_{\ \,\NCRep_n(A)} &       \\
\Alg_k    &                          &
\end{diagram}
\]
The functor $\,\NCRep_n(A) \,$ is defined by the same formula as $ \Rep_n(A) $ in \eqref{rep},
but with $ B $ being an associative algebra. It turns out that $ \NCRep_n(A) $ is representable,
and quite remarkably, its representing object $\,\sqrt[n]{A}\,$ has a simple and explicit algebraic
construction\footnote{This construction has been recently used in \cite{LBW} (see also \cite{L}),
from which we borrow the notation $\sqrt[n]{A}$.}:
\begin{equation}
\la{sqrt}
\sqrt[n]{A} = (A \ast_{k} \M(n))^{\M(n)}\ ,
\end{equation}
where $\,A \ast_{k} \mathbb{M}(n)\,$ is the free product of $ A $ and $\,\M(n) := \M(n,k) $ as $k$-algebras (i.e.,
the coproduct in the category $ \Alg_k$) and $\,(\mbox{---})^{\M(n)} $ denotes
the centralizer of $ \M(n) $ in $ A \ast_{k} \M(n)\,$. The algebra $\sqrt[n]{A}$ is thus the universal coefficient ring for the
$n$-dimensional associative representations of $ A $; it can be thought of as the coordinate ring of a noncommutative affine
scheme in the sense of \cite{C}.

Now, the inclusion functor $\,\cAlg_k \into \Alg_k \,$ has an obvious left adjoint, which is abelianization:
$\, A \mapsto A_{\natural} := A/A[A,A]A \,$, so representability of \eqref{rep} follows immediately
from \eqref{sqrt}: the commutative algebra $\, A_n := k[\Rep_{n}(A)] \,$ representing  $ \Rep_{n}(A) $ is given by the formula
\begin{equation}
\la{sqrta}
A_n = (\sqrt[n]{A}\,)_{\natural}\ .
\end{equation}

The aim of this paper is to construct the (non-abelian) derived functors of the functor $\, A \mapsto \Rep_n(A) \,$ in the sense of Quillen \cite{Q}. To do this, we first extend \eqref{rep} to the category $ \DGA_k $ of differential graded (DG) $k$-algebras, defining
$ \Rep_{n}(R) $ for a fixed $ R \in \DGA_k $ by
\begin{equation}
\la{repR}
\Rep_{n}(R):\,\cDGA_{k} \to \Sets\ ,\quad
B \mapsto \Hom_{\DGA}(R,\,B \otimes_{k} \M(n)) \ ,
\end{equation}
where $\,\cDGA_{k}\,$ is the category of commutative differential graded algebras. To represent $\Rep_{n}(R)$ we then proceed
as in the case of usual algebras: first, we prove representability in the category of all DG algebras, and then by abelianizing, we get representability of \eqref{repR} in the category of commutative DG algebras. It turns out that again the natural differential graded analogues of algebras \eqref{sqrt} and \eqref{sqrta}, to wit $\,\sqrt[n]{R} := (R \ast_{k} \M(n))^{\M(n)} \,$  and $\,R_n := (\sqrt[n]{R})_{\natural}\,$,  represent the corresponding functors explicitly.

Given an associative algebra $\,A \in \Alg_k $, we now take its almost free DG resolution $\,R \to A \,$ and define
$$
A_n^\bullet := {\rm H}^{\bullet}(R_n)\ .
$$
We prove that the assignment $\,A \mapsto A_n^\bullet \,$ is independent of the choice of resolution, and in fact, defines a functor $\,\DRep^{\bullet}_{n}\!:\,\Alg_k \to \gcAlg_k \,$ with values in the category of graded commutative algebras, so that
$\,\DRep^{0}_{n}(A) \cong \Rep_n(A)\,$ (see Section~\ref{S2.1} below).

The idea of deriving representation schemes is certainly not new: there are several different (and some, in fact, more general)
approaches in the literature: see, e.g., \cite{BCHR, T, TV}, and especially \cite{CK}.
Our construction is motivated by recent developments in noncommutative
geometry (see \cite{KR, G, L, CEG, vdB, vdB1, Be}), which are based on Kontsevich's idea \cite{K} that the family of representation
schemes $ \{\Rep_n(A)\} $ should be viewed as a `good approximation' of the noncommutative `$\Spec(A)$'. This idea turns out to be very fruitful in practice, since it allows one to find correct definitions for various structures on a noncommutative algebra $ A $ by requiring that these structures induce standard geometric structures on all representation spaces $ \Rep_n(A) $. Van den Bergh \cite{vdB} has recently proposed a concrete realization of this principle, by introducing a natural functor $\,(\mbox{---})_n:\, \Bimod(A) \to \Mod(A_n)\,$,
which transforms noncommutative objects on $A$ (viewed as bimodules over $A$) to the corresponding classical objects on $ \Rep_n(A) $.
We will construct a functor $\,\D(\mbox{---})_n:\,\D(\Bimod\,A) \to \D(\DGMod\,R_n)\,$ from the derived category of bimodules over $A$ to the derived category of DG modules over $ R_n $, where $\, R \to A\,$ is an almost free resolution of $A$.
This functor is (essentially) independent of the choice of resolution and should be viewed as a derived functor of
$\,(\mbox{---})_n\,$ in the sense of differential homological algebra (see Section~\ref{S2.2}). When combined with cohomology, it yields
a functor $\,(\mbox{---})^\bullet_n := \H^\bullet\D(\mbox{---})_n: \,\D(\Bimod\,A) \to \GrMod(A^\bullet_n)\,$, which transforms bimodules over $A$ to graded modules over $ A_n^\bullet $. Since passing from $ \Rep_n(A) $ to the DG scheme $\,\Rep_n(R) \,$ amounts (in a sense) to a desingularization of $ \Rep_n(A) $, one should expect that $\,\D(\mbox{---})_n\,$ will play a role in geometry of arbitrary (in particular, {\it homologically}\, smooth) algebras similar to the role of Van den Bergh's functor in geometry of {\it formally} smooth algebras. As a simple illustration of this idea, we compute the derived tangent spaces $ \T^\bullet_\varrho\DRep_{n}(A) $ of $\,\DRep^{\bullet}_{n}(A)\,$ at $\, \varrho \in \Rep_n(A) $, which turn out to
be isomorphic to Hochschild cohomology of the representation $ \varrho: A \to \M(n) \,$ (see Section~\ref{S2.3}):
\begin{equation}
\la{dtan}
\T^0_\varrho\DRep_{n}(A) \cong \T_\varrho \Rep_n(A) \cong  \Der_k(A,\,\M(n))\quad \mbox{and}\quad
\T^i_\varrho\DRep_{n}(A) \cong \H^{i+1}(A,\,\M(n))\ ,\ \forall\,i \ge 1 \ .
\end{equation}

Our construction of $ \DRep $ works naturally in greater generality, when instead of a single $k$-vector space, we take a complex
$V$ of finite total dimension. In the special case, when $V$ is concentrated in degree $ 0 $, our results agree with the earlier results of
Ciocan-Fontanine and Kapranov \cite{CK}: specifically, if $ R $ is a $\Z_{-}$-graded DG algebra, the affine DG action scheme $ \RAct(R,\,V) $ introduced in \cite{CK}, Section~3.3, is naturally isomorphic to our DG representation scheme $ \Rep_V(R) $. We prove this by verifying
that $ \RAct(R,\,V) $ and $ \Rep_V(R) $ satisfy the same universal property, although the a priori constructions of these schemes seem very different. The derived tangent spaces of $ \RAct(R,\,V) $ have been also computed in \cite{CK}, and the result (see {\it loc. cit.}, Proposition~3.5.4) agrees with \eqref{dtan}. Our method of computing $ \T^\bullet_{\varrho} \DRep_{V}(A) $ using Van
den Bergh's functor is different from \cite{CK} and apparently quite a bit simpler (cf. Section~\ref{S2.3} below).

The main advantage of our approach is the explicitness of algebraic constructions, which should allow concrete computations. We will give several examples in the end of the paper (see Section~\ref{S3}); these examples are chosen more or less at random, with a sole goal to illustrate the theory.
\subsection*{Acknowledgements}{\footnotesize{
The first author is very grateful to Andrei Okounkov who raised interesting questions related to \cite{BC}; attempts to clarify
his questions have been a motivation behind this work. We are also very grateful to Martin Kassabov for several insightful suggestions and
to Frank Moore, whose computer DG algebra package for {\tt Macaulay2} we have been using extensively.
We also thank Frank for his assistance with computations presented in Section~\ref{S3}.
The first author was partially supported by the NSF grant DMS 09-01570. The second author acknowledges the support by a
NSF Research Fellowship, and third author is currently funded by the Swiss National Science Foundation (Ambizione
Beitrag Nr. PZ00P2-127427/1).}}

\section{Main Results}
\la{S2}
\subsection{}
\la{S2.1}
Let $\,V = [\ldots \to V^i \to V^{i+1} \to \ldots ]\,$ be a complex of $k$-vector spaces, with $\,\sum_i \dim_k V^i < \infty \,$.
The graded endomorphism ring $ \End \, V $ is then naturally a DG algebra. Using $\End \,V $, we define the functor
\begin{equation*}
\la{root}
\sqrt[V]{-}:\ \DGA_{k} \to \DGA_{k}\ ,\quad R \mapsto \sqrt[V]{R}:=(R \ast_k \End\,V)^{\End\,V} \ ,
\end{equation*}
where $\,R \ast_k \End\,V \,$ is the coproduct in the category $\, \DGA_k \,$ and
$$
(R \ast_k \End\,V)^{\End\,V} :=
\{w \in R \ast_k \End\,V \ :\ [w,\,m]=0 \ ,\ \forall \, m \in \End\,V\}\ ,
$$
with commutators being taken in the graded sense. For any $\,R \in \DGA_k \,$, we denote by $\,R_\natural := R/R[R,R]R \,$
the abelianization of $ R $, and set
\begin{equation}
\la{rv}
R_V :=(\!\sqrt[V]{R})_\natural \ .
\end{equation}
The following lemma is a generalization  of a classic result of Cohn (see \cite{C}, Sect.~6, formula (2)).
\blemma
\la{L1}
For any $ R,\,S \in \DGA_k $ and $\,B \in \cDGA_k \,$, there are canonical isomorphisms

$(a)$ $\,\Hom_{\DGA}(\sqrt[V]{R},\,S) \cong \Hom_{\DGA}(R,\,S \otimes \End\,V) \,$,

$(b)$ $\, \Hom_{\cDGA}(R_V,\,B) \cong \Hom_{\DGA}(R,\, B \otimes \End\,V)\,$.

\elemma
\noindent
The proof of Lemma~\ref{L1}$(a)$ given in \cite{C} for ordinary algebras extends (with some more or less trivial modifications) to all DG algebras, and part $(b)$ is immediate from part $(a)$. It follows from $(b)$ that the commutative DG algebra $ R_V $ represents the functor

\begin{equation*}
\la{DGrep}
\Rep_{V}(R):\,\cDGA_{k} \to \Sets\ ,\quad
B \mapsto \Hom_{\DGA}(R,\,B \otimes \End\,V) \ ,
\end{equation*}
and thus should be thought of as the coordinate ring $ k[\Rep_V(R)] $ of an affine DG scheme $ \Rep_V(R) $ (cf. \cite{CK}, Sect.~2.2).

Now, recall that every algebra $ A \in \Alg_k $ has an {\it almost free resolution} in the category $ \DGA_k $, which is given by a quasi-isomorphism $\,\,R \to A \,$. Here, $ R $ is a DG algebra, whose
underlying graded algebra $\,|R| = \oplus_{i \in Z} R^i \,$ is free,
with $ R^i = 0 $ for all $ i > 0 $. Given two almost free resolutions
$\,R_1\to A\,$ and $\,R_2 \to B\,$ and an algebra map $\,f:\,A \to B\,$, there is a homomorphism
$\,\tilde{f}:\,R_1 \to R_2\,$ in $ \DGA_k $, such that $\,\H^0(\tilde{f})=f\,$; moreover, $\,\tilde{f}\,$
is unique up to (multiplicative) homotopy (see \cite{Q}, Ch.~I, or \cite{CK}, Sect.~3.6). Using these
(and other) standard results of homotopical algebra, one can verify that the functor $\,R \mapsto R_V \,$ defined by \eqref{rv}
preserves quasi-isomorphisms. As often happens with non-abelian derived functors, this verification requires some
technical preparations and is not immediate. The consequence is the following theorem, which is our first main result.
\begin{theorem}
\la{T1}
Let $ A $ be an associative unital $k$-algebra.

$(a)$\ For any almost free resolutions $ R_1 $ and $ R_2 $ of $A$, there is a
quasi-isomorphism $ f: (R_1)_V \to (R_2)_V $.

$(b)$\ The assignment $\,A \to A_V^\bullet := \H^{\bullet}(R_V)\,$ defines a functor $\,\Alg_{k} \to \gcAlg_k\,$, which
is independent of the choice of almost free resolution $\,R \to A\,$.

$(c)$ If $V$ is concentrated in degree $0$, then $\,A_V^0 \cong A_V \,$.
\end{theorem}
Since $\,\Spec(A_V) = \Rep_V(A) \,$ when $V$ is concentrated in degree $0$, Theorem~\ref{T1}$(c)$ implies
that the graded scheme $ \boldsymbol{\Spec}(A_V^\bullet)\,$ should be viewed as a derived representation scheme
of $A$. To simplify the notation, we set $\, \DRep^\bullet_{V}(A) := \boldsymbol{\Spec}(A_V^\bullet) $ and write
$ \DRep^\bullet_{n}(A) $ instead of $ \DRep^\bullet_{V}(A) $ when $\,V= k^n\,$.

\subsection{}
\la{S2.2}
From now on, we assume that $ V $ is concentrated in degree $0$. For an algebra $ A \in \Alg_k $, we write
$\, A_V := k[\Rep_V(A)]\,$ and  let
$\,\pi_V: A \to A_V \otimes \End\,V\,$
denote the universal representation. Restricting scalars via $\pi_V$, we can regard $\,A_V \otimes \End\,V\,$ as a
bimodule over $ A $, or equivalently, as a left module over the enveloping algebra $\,\eA := A \otimes A^{\rm opp} \,$.
Since $ A_V $ is commutative, the image of $ A_V $ under the natural inclusion $\,A_V \into A_V \otimes \End\,V \,$
lies in the center of this bimodule. Hence, we can regard $\,A_V \otimes \End\,V\,$ as $\,\eA$-$A_V$-bimodule. Now, following Van den Bergh (see \cite{vdB}, Sect.~3.3), we define the additive functor
\begin{equation}
\la{E1}
(\mbox{---})_V:\ \Bimod(A) \to \Mod(A_V)\ , \quad M \mapsto M \otimes_{\eA} (A_V \otimes \End\,V)\ .
\end{equation}
As mentioned in the Introduction, this functor plays a key role in noncommutative geometry of smooth algebras, transforming noncommutative objects on $A$ to the classical geometric objects on $ \Rep_V(A) $. Our aim is to construct the higher derived functors of \eqref{E1}, which should replace \eqref{E1} when $A$ is not smooth. We begin by extending the Van den Bergh functor to the world of DGAs.

Fix $\, R \in \DGA_k \,$ and let $\,\pi:\,R \to \rtv{R} \otimes \End\,V\,$ denote the universal DG algebra homomorphism
corresponding to the identify functor on $ \rtv{R} $, see Lemma~\ref{L1}$(a)$. The complex
$\,\rtv{R} \otimes V \,$ is naturally a left DG module over $ \rtv{R} \otimes \End\,V $ and right DG module
over $ \rtv{R} \,$, so restricting the left action via $ \pi $ we can regard  $\,\rtv{R} \otimes V \,$ as DG bimodule over $R$ and $\rtv{R}$.
Similarly, we can make $\,V^* \otimes \rtv{R}$ a $\,\rtv{R}$-$R$-bimodule.
Using these bimodules, we define the functor
\begin{equation}
\la{rtvm}
\rtv{-}:\ \DGBimod(R) \to \DGBimod(\rtv{R}) \ , \quad
M \mapsto (V^* \otimes \sqrt[V]{R}) \otimes_R M \otimes_R (\sqrt[V]{R} \otimes V) \ .
\end{equation}
Now, recall that $\,R_V := (\rtv{R})_\natural\,$ is a commutative DGA. Using the natural projection $\,\rtv{R} \to R_V \,$, we regard
$R_V$ as a DG bimodule over $ \rtv{R} $ and define
\begin{equation}
\la{abb}
(\mbox{---})_\natural:\ \DGBimod(\rtv{R}) \to \DGMod(R_V) \ , \quad
M \mapsto M_\natural := M \otimes_{(\rtv{R})^{\rm e}} R_V  \ ,
\end{equation}
which is nothing but the abelianization functor on bimodules. Combining
\eqref{rtvm} and \eqref{abb}, we define
\begin{equation}
\la{E11}
(\mbox{---})_V:\ \DGBimod(R) \to \DGMod(R_V)\ , \quad M \mapsto
M_V := (\rtv{M})_\natural \ .
\end{equation}
As suggested by its notation, the functor \eqref{E11} is a DG extension of \eqref{E1}. In fact, if $ R = A $ is
a DG algebra with a single nonzero component in degree $0$, the category $ \Bimod(A) $ can be viewed as
a full subcategory of $\DGBimod(R)$ consisting of bimodules concentrated in degree $0$.
It is easy to check then that the restriction of \eqref{E11} to this subcategory coincides with \eqref{E1}.

The next lemma is analogous to Lemma~\ref{L1} for DG algebras; it holds, however, in greater generality:
for homomorphism complexes $\, \HOM^\bullet $ of DG modules. We recall that, if $ R $ is a DG algebra and
$M$,$\,N$ are DG modules over $R$, $\, \HOM^\bullet_R(M,\,N) \,$ is a complex of vector spaces with
$n$-th graded component consisting of all $R$-linear maps $\,f:\, M \to N\,$ of degree $ n $ and
the $n$-th differential given by $\,d(f) = d_N \circ f - (-1)^n f \circ d_M\,$.

\blemma
\la{L2}
There are canonical isomorphisms of complexes

\vspace{0.8ex}

$(a)$ $\,\HOM^\bullet_{(\rtv{R})^{\rm e}}(\rtv{M},\,N)
\cong \HOM^\bullet_{\eR}(M,\,N \otimes \End\,V) \,$,

\vspace{0.8ex}

$(b)$ $\, \HOM^\bullet_{R_V}(M_V,\,N) \cong \HOM^\bullet_{\eR}(M,\,N \otimes \End\,V) \,$.

\elemma

\begin{example}
Let $ R $ be a DG algebra. Denote by $ \Omega^{1} R $ the
kernel of the multiplication map $\,R \otimes R \to R\,$.
This is naturally a DG bimodule over $ R $, which, as in the case
of ordinary algebras, represents the complex of $k$-linear graded
derivations $\,\DER_k^\bullet(R,\,M) \,$ (see, e.g., \cite{Q1}, Sect.~3.1 and~3.2).
Thus, for any $\, M \in \DGBimod(R) \,$, there is a canonical isomorphism of complexes of vector spaces
\begin{equation}
\la{der}
\DER_k^\bullet(R,\,M) \cong \HOM^\bullet_{\eR}(\Omega^{1} R,\,M)\ .
\end{equation}
Using \eqref{der} and Lemma~\ref{L2}, one can establish canonical isomorphisms %
\begin{equation}
\la{omvdb}
\rtv{\Omega^1 R} \cong \Omega^1(\rtv{R}) \quad , \quad
(\Omega^1 R)_V \cong  \Omega^1(R_V)\ .
 \end{equation}
This should be compared to \cite{vdB}, Proposition~3.3.4.
\end{example}

Next, we recall that if $ R $ is a DG algebra, every
DG module $M$ over $ R $ has a {\it semi-free resolution}
$\,F \to M \,$, which is similar to a free (or projective)
resolution for ordinary modules over ordinary algebras
(see \cite{FHT}, Sect.~2). To construct the derived
functors of \eqref{E1} we now follow the standard procedure
in differential homological algebra.

Given an algebra $ A \in \Alg_k $ and a complex $M$ of bimodules over $A$, we first
pick an almost free resolution $\,f: R \to A\,$ in $\, \DGA_k \,$ and consider $ M $
as a DG bimodule over $ R $ via $f$. Then, we pick a semi-free resolution
$\,F(R,M) \to M \,$  in the category $ \DGBimod(R) $ and apply to $F(R,M)$
the functor \eqref{E11}. The result is described by the following theorem,
which is the second main result of this paper.
\begin{theorem}
\la{T2}
Let $ A $ be an associative $k$-algebra, and let $ M $ be a complex of bimodules over $A$.

$(a)$\ The assignment $\,M \mapsto F(R,M)_V\,$ induces a well-defined functor between the derived categories
$$
\D(\mbox{---})_V:\ \D(\Bimod\,A) \to \D(\DGMod\,R_V)\ ,
$$
which is independent of the choice of the resolutions $ R \to A $ and $ F \to M $ up to auto-equivalence of 
$\, \D(\DGMod\,R_V) $ inducing the identity on cohomology.

$(b)$\ Taking cohomology $\,M \mapsto \H^\bullet\D(M)_V\,$ yields a functor
$$
(\mbox{---})^\bullet_V:\ \D(\Bimod\,A) \to \GrMod(A_V^\bullet)\ ,
$$
which depends only on the algebra $A$ and the vector space $ V $.

$(c)$\ If $ M \in \Bimod(A) $ is viewed as a 0-complex in $\D(\Bimod\,A)$, then $\,M_V^0 \cong M_V \,$.
\end{theorem}
The proof of part $(a)$ is standard differential homological
algebra; the fact that $ \D(\mbox{---})_V $ is independent of
resolutions follows from results of Keller (see, e.g.,
\cite{Ke1}, Sect.~8.4). Part $(b)$ is immediate from $(a)$,
and $(c)$ is proved by direct computation.

\subsection{}
\la{S2.3}
We will use the above construction to compute the derived tangent spaces $ \T^\bullet_\varrho \DRep_{V}(A) $ at $\,\varrho \in \Rep_{V}(A)\,$.
It is instructive to compare our computation with the ones in \cite{vdB}, see {\it loc. cit.}, Section~3.3.

Recall that if $X$ is an ordinary affine $k$-scheme and
$ x \in X $ is a $k$-point, the tangent space
$ \T_x X $ to $ X $ at $ x $ is defined by
$\,\T_x X := \Der_k(k[X],\,k_x)\,$.  Similarly (cf. \cite{CK}, (2.5.6)), if $ \X $ is a DG scheme, and $ x \in \X^0 $ is a $k$-point, the
tangent DG space to $ \X $ at $x$ is defined by
$$
\T_x\,\X := \DER_k^\bullet(k[\X],\,k_x)\ .
$$
Now, if $ X^\bullet := \boldsymbol{\Spec}\,\H^\bullet(k[\X]) $ is
the underlying derived scheme of $ \X $, the {\it derived} tangent space $ \T^\bullet_x X $ to $ X^\bullet $ is given, by definition,
by cohomology:
\begin{equation}
\la{tspace}
\T^\bullet_x X := \H^\bullet(\DER_k(k[\X],\,k_x))\ .
\end{equation}

Let $\,\rho:\,R_V \to k \,$ be the DG algebra
homomorphism corresponding to a $k$-point in $ \DRep_V^\bullet(A) $, and let $\,\varrho:\,R \to \End(V)\,$ be the representation
corresponding to $ \rho $. Then, we have canonical isomorphisms of complexes of vector spaces
\begin{eqnarray*}
\DER_k^\bullet(R_V, \, k) &\cong& \HOM^\bullet_{R_V}(\Omega^1(R_V),\,k) \\
                 &\cong& \HOM^\bullet_{R_V}((\Omega^1 R)_V,\,k)\qquad [\mbox{see}\,\eqref{omvdb}] \\
                 &\cong& \HOM^\bullet_{R_V} (\rtv{\Omega^1(R)} \otimes_{(\rtv{R})^{\rm e}} R_V,\,k) \\
                 &\cong& \HOM^\bullet_{(\rtv{R})^{\rm e}} (\rtv{\Omega^1(R)},\ \HOM^\bullet_{R_V}(R_V,\,k)) \\
                 &\cong& \HOM^\bullet_{(\rtv{R})^{\rm e}} (\rtv{\Omega^1(R)},\,k) \\
                 &\cong& \HOM^\bullet_{\eR} (\Omega^1 R ,\,\End\,V) \qquad [\mbox{see\,Lemma~\ref{L2}$(a)$}]\\
                 &\cong& \DER_k^\bullet(R,\,\End\,V)\ ,
\end{eqnarray*}
which imply
$$
\T^\bullet_\varrho\DRep_{V}(A) := \H^\bullet[\DER_k(R_V, \, k)]
\cong \H^\bullet[\DER_k(R,\,\End\,V)]\ .
$$
The following proposition is now an immediate consequence of \cite{BP}, Lemma~4.2.1 and Lemma~4.3.2.
\bprop[cf. \cite{CK}]
\la{P2}
$$
\T^i_\varrho\DRep_{V}(A) \cong \left\{
\begin{array}{lll}
\Der_k(A,\,\End\,V)\ & \mbox{\rm if} &\ i = 0\\*[1ex]
\H^{i+1}(A,\,\End\,V)\ & \mbox{\rm if} &\ i \ge 1
\end{array}
\right.
$$
\eprop

\begin{remark}
As mentioned in the Introduction, in case when $V$ is a single vector space and $R \in \DGA_k$ is almost free,
one can show that $ \Rep_V(R) $ is isomorphic the DG scheme $ \RAct(R,\,V) $ constructed in \cite{CK}.
This implies that $\,\T^i_\varrho\DRep_{V}(A)\,$ should be isomorphic to $\,\T^i_\varrho \RAct(R,\,V)\,$,
which is indeed the case, as one can easily see by comparing our Proposition~\ref{P2} to \cite{CK},
Proposition~3.5.4$(b)$.
\end{remark}

\section{Examples}
\la{S3}
In this section, we assume that $V=k^n$. Given an explicit almost free resolution $R \to A $, the DG algebra $\,
R_V = (\rtv{R})_\natural $ can be described explicitly. Specifically, let $\{r_i\}_{i \in I}$ be a
set of generators of $R$, and let $d_R$ be its differential. Consider
a free graded algebra $\tilde{R}$ with generators
$\,\{r_{i}^{jk}\}_{i\in I,1\leq j,k\leq n}\,$ of degree $\,|r_{i}^{jk}|=|r_{i}|\,$. Define linear functions
$\,\{f_{jk}\}_{1\leq j,k\leq n}\,$ on (tensor) products of generators of $R$ by
$$
f_{jk}:\ R \to \tilde{R}\quad , \qquad
r_{t_{1}}r_{t_{2}}...r_{t_{m}}\ \mapsto\ \sum\ r_{t_{1}}^{js_{1}}\,
r_{t_{2}}^{s_{1}s_{2}}\, \ldots\, r_{t_{m}}^{s_{m-1}k}\ ,
$$
where the sum is taken over all $\, 1\leq s_{1},\,\ldots\,,s_{m-1}\leq n\,$,
and extend this by linearity to the whole of $R$. Now, define a differential $d$ on generators of $\tilde{R}$ by
$$
d(r_{i}^{jk}) := f_{jk}(d_{R}(r_{i}))\ .
$$
and extend it the Leibniz rule to whole of $\tilde{R}$. This makes $\,\tilde{R}\,$ a DG algebra. The abelianization
of $\tilde{R}$ is then a free (graded) commutative algebra with generators $\,r_i^{jk}\,$ and the differential induced
by the differential of $\tilde{R}$. Then, we have the following result.
\bthm
\la{comp}
There is an isomorphism of DG algebras $\, \tilde{R} \cong \sqrt[V]{R} \,$. Consequently, $\,R_V \cong \tilde{R}_\natural\,$.
\ethm

Using Theorem~\ref{comp}, we can construct explicitly a finite presentation of the graded algebra
$\, A_V^\bullet \,$, whenever we have a finite almost free resolution of $\, A\,$.
In practice, for many interesting algebras, such resolutions are available. For example, the DG algebras introduced
recently in \cite{G1} and \cite{Ke} provide very interesting (in a sense, canonical) finite resolutions for all
3D Calabi-Yau algebras, including $\,k[x, y, z]\,$, $\,{\mathcal U}(\mathfrak{sl}_{2})\,$, Sklyanin algebras, and many others.
This allows us to describe the corresponding derived representation schemes quite explicitly, which is rather unusual for
derived functors in homotopical algebra.

\begin{example}
Let $ \mathcal{U}(\mathfrak{sl}_{2}) $ be the universal
enveloping algebra of the Lie algebra $ \mathfrak{sl}_{2}(k) $. By \cite{G1}, Example~1.3.6, it has
a finite DG resolution $\,R \to \mathcal{U}(\mathfrak{sl}_{2})\, $, where $ R $ is the free graded algebra
$$
R = k \langle x,\,y,\,z;\, X,\, Y,\, Z;\, t\rangle\ ,
$$
with generators $\,x,\,y,\,z\,$ having degree $0$;  $\,X,\,Y,\,Z\,$ having degree $-1$ and $t$ having
degree $-2$. The differential on $R $ is defined by
$$
dX = yz-zy + x\ ,\quad dY = zx-xz + y\ ,\quad dZ = xy-yx + z\ ,\quad dt = [x, X] + [y, Y] + [z, Z]\ .
$$
Theorem~\ref{comp} then implies that
$$
R_n = k[x_{ij}\,,\,y_{ij}\,,\,z_{ij}\,;\,X_{ij}\,,\,Y_{ij}\,,\,Z_{ij}\,;\,t_{ij} \ |\ 1 \leq i,\,j \leq n]\ ,
$$
where the generators $\,x_{ij},\,y_{ij},\,z_{ij}\,$ have degree zero,
$\,X_{ij},\,Y_{ij},\,Z_{ij}\,$ have degree $ -1$, and $\,|t_{ij}| = -2\,$. The differential on $\,R_n\,$ is given by
\begin{eqnarray*}
dx_{ij}&=&dy_{ij}=dz_{ij}=0\ ,\\
dX_{ij}&=& \sum_{k=1}^n (y_{ik}z_{kj}-z_{ik}y_{kj})+ x_{ij}\ ,\\
dY_{ij}&=& \sum_{k=1}^n (z_{ik}x_{kj}-x_{ik}z_{kj})+ y_{ij}\ ,\\
dZ_{ij}&=& \sum_{k=1}^n (x_{ik}y_{kj}-y_{ik}x_{kj})+ z_{ij}\ ,\\
dt_{ij}&=& \sum_{k=1}^n x_{ik}X_{kj}-X_{ik}x_{kj}+y_{ik}Y_{kj}-Y_{ik}y_{kj}+
z_{ik}Z_{kj}-Z_{ik}z_{kj}\ .
\end{eqnarray*}
\end{example}

We now give a simple example when one can actually compute the cohomology of $\, R_n\,$.
\begin{example}
Let $\,A = k[x,y]\,$ be the (commutative) polynomial algebra in two variables. One can easily check that $A$ has resolution
$\,R := k\langle x,\,y;\,t\rangle\,$, with generators $\,x,y\,$ in degree $0$, $\,t$ in degree $-1$, and $\,dt=xy-yx\,$.
In this case, using Gr\"{o}bner basis techniques, we can compute
$$
k[x,y]^\bullet_2  = \frac{k[x,y]_2\,\langle r,\,s \rangle}{\left(\begin{array}{c}
x_{21}r-y_{21}s,
x_{12}r-y_{12}s,\\
(x_{11}-x_{22})r-(y_{11}-y_{22})s,\\
sr, \; rs,\; r^{2},\; s^{2}\end{array}\right)}\ ,
$$
where the two extra generators $r$ and $s$ have degree $-1$.
\end{example}
Similar presentations can be constructed for other polynomial algebras, using a computer algebra software.
For example, $\, k[x,y]^\bullet_3\,$ has $6$ generators in nonzero degrees, all of them in degree $-1$; it has nonzero components
only in degrees $\,0,\, -1,\, -2,\,$ and $-3$. The algebra $\,k[x,y,z]^\bullet_2\,$ has $16$ generators in degree
$-1$, and others in lower degrees; regarded as a $\, k[x,y,z]_2$-module, its minimal generating set has
$16$ elements in degree $-1$; $\, 56$ in degree $-2$; $\,128 $ in degree $-3$;\,
$233$ in degree $-4$, and more in lower degrees.

\end{document}